\documentclass[amsmath,amssymb,10pt]{article}
\usepackage{amsmath, amsfonts, amssymb, amsthm}

\begin{document}

%
\setcounter{page}{1}  
%

\newcommand{\N}{\mathbb{N}}
\newcommand{\Z}{\mathbb{Z}}

\newcommand{\R}{\mathbb{R}}

\newcommand{\capnab}{\nabla^{\nu}_{a^*}}

 \newcommand{\xproof}[1]{\par\noindent{\em Proof.  }\ #1 \, \nobreak \hfill $\blacksquare$ \vspace{.3pc}}
 \newtheorem{definition}[theorem1]{Definition}  
 \newtheorem{example}[theorem1]{Example}
\newtheorem{remark}[theorem1]{Remark}
\newtheorem{claim}[theorem1]{Claim}
 \newtheorem{lemma}[theorem1]{Lemma}  
 \newtheorem{theorem}[theorem1]{Theorem} 
 \newtheorem{corollary}[theorem1]{Corollary} 

\begin{center}

\vspace{0.4cm} {\bf{\Large  Initial 
and Boundary Value Problems for the Caputo Fractional Self-Adjoint Difference 
Equations}}
\\
\vspace{0.3cm}
Kevin Ahrendt \\
University of Nebraska-Lincoln \\
Department of Mathematics \\
Lincoln, NE 68588-0130 USA \\
{\tt kahrendt@huskers.unl.edu} \\

\vspace{0.3cm}

Lydia DeWolf\\
Department of Mathematics\\
Union College\\
Jackson, TN, 38305 USA \\
{\tt  lydia.dewolf@my.uu.edu}
\vspace{0.3cm}

Liam Mazurowski\\
Department of Mathematical Sciences\\
Carnegie Mellon University\\
Pittsburgh, PA 15213 USA \\
{\tt lmazurow@andrew.cmu.edu}
\vspace{0.3cm}

Kelsey Mitchell\\
Department of Mathematics\\
 Buena Vista University \\
Storm Lake, IA 50588 USA \\
{\tt mitckel@bvu.edu}
\vspace{0.3cm}

Tim Rolling \\
Department of Mathematics \\
University of Nebraska-Lincolnt\\
Lincoln, NE 68588-0130 USA \\
{\tt trolling2@unl.edu} \\
\vspace{0.3cm}

Dominic Veconi\\
Department of Mathematics \\
Hamilton College \\
Clinton, NY 13323 USA \\
{\tt dveconi@hamilton.edu}

\vspace{0.4cm} {\it Communicated by Allan Peterson} 
\end {center}
\newpage
\begin{abstract}
In this paper we develop the theory of initial and boundary value problems for 
the self-adjoint nabla fractional difference equation containing a Caputo 
fractional nabla difference that is given by
\[
\nabla[p(t+1)\nabla_{a*}^\nu x(t+1)] + q(t)x(t) = h(t),
\] where $0 < \nu \leq 1$. We give an introduction to the nabla 
fractional calculus with Caputo fractional differences. We investigate 
properties of the specific self-adjoint nabla fractional difference equation 
given above. We prove existence and uniqueness theorems for both initial and 
boundary value problems under appropriate conditions. We introduce the 
definition of a Cauchy function which allows us to give a variation of constants 
formula  for solving  initial value problems. We then show that this Cauchy 
function is important in finding  a Green's function for a boundary value 
problem with Sturm-Liouville type boundary conditions. Several inequalities 
concerning a certain Green's function are derived. These results are important 
in using fixed point theorems for proving the existence of solutions to boundary 
value  problems for nonlinear fractional equations related to our 
linear self-adjoint equation.\\


\end{abstract}

\section{Preliminary definitions and theorems}
In this section we will develop the basic definitions and theorems needed for our results. For an introduction to whole order difference calculus, see Kelley and Peterson \cite{Kelleydiff}. For more information in the fractional results using the backwards difference operator, see Hein et al. \cite{Hein} and Ahrendt et al. \cite{ahrendt}.

\begin{definition}\label{defna} \cite{ahrendt}
Let $a\in \R$.  Then the set $\N_a$ is given by $\{a,a+1,a+2,\hdots\}$. Furthermore, if $b \in \N_a$, then $\N_a^b$ is given by $\{a,a+1, \ldots, b-1,b\}$.
\end{definition}

\begin{definition}\label{defnabdiff} \cite{ahrendt}
Let $f:\N_a \to \R$.  Then the \textit{nabla difference of $f$} is defined by 
\[
(\nabla f)(t) := f(t) - f(t-1),
\]
for $t\in \N_{a+1}$.  For convenience, we will use the notation $\nabla f(t) := (\nabla f)(t)$. For $N \in \N$, we have that the \textit{$N$th order fractional difference} is recursively defined as
\[
\nabla^N f(t) := \nabla (\nabla^{N-1} f(t)),
\]
for $t\in \N_{a+N}$.  
\end{definition}

\begin{definition} \label{defint} \cite{ahrendt}
Let $f:\N_a \to \R$ and let $c,d\in \N_a$.  Then the \textit{definite nabla integral of $f$ from $c$ to $d$} is defined by
\[
\int_c^d f(s)\nabla s := \begin{cases}
					\sum_{s=c+1}^d f(s), & c<d, \\
					0, & d\leq c.
				\end{cases}
\]
\end{definition}

\begin{definition} \label{defrising} \cite{ahrendt}
Let $t\in \R$ and let $n\in\Z^+$.  Then the \textit{rising function} is defined by
\[
t^{\overline n} := (t)(t+1)\cdots (t+n-1) = \frac{\Gamma(t+n)}{\Gamma(t)},
\]
where $\Gamma$ is the gamma function. For $\nu \in \R$, the \textit{generalized rising function} is then defined by
\[
t^{\overline \nu} := \frac{\Gamma(t+\nu)}{\Gamma(t)}
\]
for $t$ and $\nu$ such that $t+\nu \not \in \{\ldots,-2,-1,0\}$. If $t$ is a non-positive integer and $t + \nu$ is not a non-positive integer then we take by convention $t^{\overline \nu} = 0$.  
\end{definition}

\begin{theorem}[Fundamental Theorem of Nabla Calculus]\label{fund} \cite{ahrendt}
Assume the function $f:\N_a^b\rightarrow\R$ and let $F$ be a nabla antidifference of $f$ on $\N_a^b$, then
	\[
	\int_a^b f(t)\nabla t= F(t)\bigg|_a^b := F(b)-F(a).
	\]
\end{theorem}

The following definitions extend the nabla difference and nabla integral to fractional value orders.

\begin{definition}\label{fracsum} \cite{ahrendt}
Let $f:\mathbb{N}_{a+1}\rightarrow\mathbb{R}$, $\nu>0$, $\nu\in\mathbb{R}$. The $\nu^{th}$ \textit{order nabla fractional sum of f} is defined as 
	\[
	\nabla_a^{-\nu}f(t):=\sum_{s=a+1}^t\frac{(t-s+1)^{\overline{\nu-1}}}{\Gamma(\nu)}f(s),
	\]
for $t\in\mathbb{N}_{a}$.
\end{definition}

\begin{definition}\cite{ahrendt}
Let $f:\N_{a+1} \to \R$, $\nu \in \R$, $\nu >0$, and $N = \lceil\nu\rceil$. Then $\nu^{th}$ \textit{order nabla fractional difference of f} is defined as
\[
\nabla_a^{\nu} f(t) := \nabla^N \nabla_a^{-N-\nu}f(t),
\]
for $t \in \N_{a+N}$.
\end{definition}

The next theorem gives results for composing nabla fractional sums and differences in certain cases.

\begin{theorem}[Composition Rules] \label{comprules} \cite{ahrendt}
Let $\mu, \nu > 0$ and let $f:\N_a \to \R$.  Set $N = \lceil \nu \rceil$.  Then
\[
\nabla_a^{-\nu}\nabla_a^{-\mu} f (t) = \nabla_a^{-\nu-\mu}f(t), \quad t\in \N_a,
\]  
and
\[
\nabla_a^\nu \nabla_a^{-\mu} f(t) = \nabla_a^{\nu-\mu} f(t), \quad t\in \N_{a+N}.
\]
\end{theorem}

While we gave the traditional definition of a nabla fractional difference above, we will focus on the Caputo nabla fractional difference for the rest of the paper and only appeal to the previous definition when needed. The following definition has been adapted from Anastassiou in \cite{anas}.
\begin{definition}\label{capdiff}
Let $f:\mathbb{N}_{a-N+1}\rightarrow\mathbb{R}$, $\nu>0$, $\nu\in\mathbb{R}$, $N=\lceil\nu\rceil$.  The $\nu^{th}$ \textit{order Caputo nabla fractional difference} is defined as 
	\[
	\nabla_{a*}^{\nu}f(t):=\nabla_a^{-(N-\nu)}(\nabla^Nf(t)),
	\]
for $t\in\mathbb{N}_{a}$. Note that the Caputo difference operator is a linear operator.
\end{definition}

\begin{theorem} [Discrete Whole-Order Taylor's Formula]\label{anasthm} \cite{anderson}
Fix $N\in \N_1$ and let $f:\N_{a-N+1} \to \R$.  Then
\[
f(t) = \sum_{k=0}^{N-1} \nabla^k f(a)\frac{(t-a)^{\overline k}}{k!} + \sum_{s=a+1}^t \frac{(t-s+1)^{\overline {N-1}}}{(N-1)!}\nabla^{N} f(s),
\]
for $t\in \N_a$.
\end{theorem}

The following theorem is adapted from Anastassiou in \cite{anas} where we use our definition of the Caputo nabla fractional difference.

\begin{theorem}[Caputo Discrete Taylor's Theorem] \label{captay} \cite{anas}
Let $\nu\in\mathbb{R}$, $\nu>0$, $N=\lceil\nu\rceil$, and $f:\N_{a-N+1}\rightarrow\R$. Then for all $t\in\N_{a}$, the representation holds
 	\[
	f(t)=\sum_{k=0}^{N-1}\frac{(t-a)^{\overline{k}}}{k!}\nabla^kf(a)+\frac{1}{\Gamma(\nu)}\sum_{\tau=a+1}^t(t-\tau+1)^{\overline{\nu-1}}\nabla_{a*}^{\nu}f(\tau).
	\]
\end{theorem}

\xproof{
Using Theorem \ref{comprules}, notice that for $f: \N_{a-N+1} \rightarrow \R$,
	\begin{align*}
	\nabla_a^{-\nu}\nabla_{a*}^{\nu} f(t) &= \nabla_a^{-\nu}\nabla_a^{-(N-\nu)}\nabla^N f(t)\\
	&=\nabla_a^{-\nu-N+\nu}\nabla^N f(t)\\
	&=\nabla_a^{-N}\nabla^N f(t), 
	\end{align*}
for $t\in\N_a$.
By Definition \ref{fracsum}, 
	\[
	\nabla_a^{-N}\nabla^N f(t)= \frac{1}{(N-1)!}\int_a^t (t-s+1)^{\overline{N-1}}\nabla^N f(s) \nabla s,
	\]
Similarly, 
	\[
	\nabla_a^{-\nu}\nabla_{a*}^{\nu} f(t)= \frac{1}{\Gamma(\nu)}\int_a^t (t-s+1)^{\overline{\nu-1}}\nabla_{a*}^{\nu} f(s) \nabla s.
	\]
So from Theorem \ref{anasthm}, 
	\begin{align*}
	f(t)&=\sum_{k=0}^{N-1} \frac{(t-a)^{\overline{k}}}{k!}\nabla^k f(a)+\frac{1}{(N-1)!}\sum_{s=a+1}^t (t-s+1)^{\overline{N-1}}\nabla^{N} f(s)\\
	&=\sum_{k=0}^{N-1} \frac{(t-a)^{\overline{k}}}{k!}\nabla^k f(a)+\frac{1}{\Gamma(\nu)}\sum_{s=a+1}^t (t-s+1)^{\overline{\nu-1}}\nabla_{a*}^{\nu} f(s),
	\end{align*}
for $t \in \N_a$, proving the result.
}

\section{Nabla fractional initial value problems}

We are interested in solutions to the nabla fractional initial value problem (IVP)
	\begin{equation} \label{eq:ckivp}
	\begin{cases}
	&\nabla_{a*}^{\nu}f(t)=h(t), \quad t\in \N_{a+1},\\ 
	&\nabla^kf(a)=c_k, \quad 0\leq k\leq N-1,
	\end{cases} 
	\end{equation}
where $a,\nu\in\R$, $\nu>0$, $N=\lceil\nu\rceil$, $c_k\in\R$ for $0\le k\le N-1$, and $f:\N_{a-N+1}\rightarrow\R$.

\begin{theorem}\label{uniq}
The solution to the IVP \eqref{eq:ckivp} is uniquely determined by 
	\[
	f(t)=\sum_{k=0}^{N-1}\frac{(t-a)^{\overline{k}}}{k!}c_k+\frac{1}{\Gamma(\nu)}\sum_{\tau=a+1}^t(t-\tau+1)^{\overline{\nu-1}}h(\tau),
	\]
for $t\in\N_{a+1}$.
\end{theorem}

\xproof{Let $f:\mathbb{N}_{a-N+1}\rightarrow\mathbb{R}$ satisfy 
	\[
	\nabla^kf(a)=c_k,
	\]
for $0\le k\le N-1$. Note that this uniquely determines the value of $f(t)$ for $a-N+1\leq t\leq a$. For $t\in\N_{a+1}$, let $f(t)$ satisfy
	\[
	\nabla^Nf(t)=h(t)-\sum_{s=a+1}^{t-1}\frac{(t-s+1)^{\overline{N-\nu-1}}}{\Gamma(N-\nu)}\nabla^Nf(s).
	\]
This recursive definition uniquely determines $f(t+1)$ from the values of $f(a-N+1),..., f(t)$, so the function is uniquely defined for all $t\in\mathbb{N}_{a-N+1}$. So for any $t\in\mathbb{N}_{a+1}$, 
	\[
	\nabla^Nf(t)+\sum_{s=a+1}^{t-1}\frac{(t-s+1)^{\overline{N-\nu-1}}}{\Gamma(N-\nu)}\nabla^Nf(s)=h(t).
	\]
Equivalently, 
	\begin{align*}
	&\frac{\Gamma(\nu)}{\Gamma(1)\Gamma(\nu)}\nabla^Nf(t)+\sum_{s=a+1}^{t-1}\frac{(t-s+1)^{\overline{N-\nu-1}}}{\Gamma(N-\nu)}\nabla^Nf(s)\\
	&=\frac{(1)^{\overline{\nu-1}}}{\Gamma(\nu)}\nabla^Nf(t)+\sum_{s=a+1}^{t-1}\frac{(t-s+1)^{\overline{N-\nu-1}}}{\Gamma(N-\nu)}\nabla^Nf(s)\\
	&=\sum_{s=a+1}^t\frac{(t-s+1)^{\overline{N-\nu-1}}}{\Gamma(N-\nu)}\nabla^Nf(s)\\
	&=\nabla_a^{-(N-\nu)}\nabla^Nf(t)\\
	&=\nabla_{a*}^{\nu}f(t)\\
	&=h(t).
	\end{align*}
Therefore, $f(t)$ solves the IVP \eqref{eq:ckivp}. Conversely, if we suppose that there is a function $f:\mathbb{N}_{a-N+1}\rightarrow\mathbb{R}$ that satisfies the IVP, reversing the above algebraic steps would lead to the same recursive definition. Therefore the solution to the IVP is uniquely defined. By the Caputo Discrete Taylor's Theorem, $f(t)$ must satisfy the representation
	\begin{align*}
	f(t)&=\sum_{k=0}^{N-1}\frac{(t-a)^{\overline{k}}}{k!}\nabla^kf(a)+\frac{1}{\Gamma(\nu)}\sum_{\tau=a+1}^t(t-\tau+1)^{\overline{\nu-1}}\nabla_{a*}^{\nu}f(\tau)\\
	&=\sum_{k=0}^{N-1}\frac{(t-a)^{\overline{k}}}{k!}c_k+\frac{1}{\Gamma(\nu)}\sum_{\tau=a+1}^t(t-\tau+1)^{\overline{\nu-1}}h(\tau),
	\end{align*}
for $t \in \N_{a+1}$.
}

\begin{example}
Solve the IVP 
	\[
	\begin{cases}
	&\nabla_{0*}^{0.7} f(t)=t, \quad t \in \N_{1},\\
	&f(0)=2.
	\end{cases}
	\]
Applying the variation of constants formula yields the following expression for $f(t)$.
	\begin{align*}
	f(t)&= \sum_{k=0}^0 \frac{t^{\overline{k}}}{k!}2 + \nabla_0^{-0.7}h(t)\\
	&=\sum_{s=1}^t \frac{(t-s+1)^{\overline{-0.3}}}{\Gamma(0.7)}s^{\overline{1}} +2\\
	\end{align*}
After summing by parts, we have
	\[f(t)=\frac{t^{\overline{1.7}}}{\Gamma\left(2.7\right)}+2.\]
\end{example}

\section{General properties of the fractional self-adjoint nabla difference equation}

For development of these properties in the continuous setting, see Kelley and Peterson \cite{Kelley}.

Let $\mathcal{D}_a:=\{x:\mathbb{N}_a\rightarrow\mathbb{R}\}$ and let the self-adjoint fractional operator $L_a$ be defined by
	\begin{equation*}
	(L_a x)(t):=\nabla[p(t+1)\nabla_{a*}^{\nu}x(t+1)]+q(t)x(t), \quad t \in N_{a+1},
	\end{equation*}
where $x\in\mathcal{D}_a, 0<\nu<1$, $p:\N_{a+1} \to (0,\infty)$ and $q:\N_{a+1}\rightarrow\R$. Note that while the operator is for values of $t \in \N_{a+1}$, the function $x$ is defined on $\N_a$. Note that $L_a$ is a linear operator.



\begin{theorem} [Existence and Uniqueness for Self-Adjoint IVPs] \label{exuq}
Let $A,B \in\mathbb{R}$, and let $h:\mathbb{N}_{a+1}\rightarrow\mathbb{R}$. The IVP
	\[
	\begin{cases}
	&L_ax(t)=h(t), \quad t \in \N_{a+1},\\
	&x(a)=A, \\
	& \nabla x(a+1)=B,
	\end{cases}
	\]
has a unique solution $x:\N_a \to \R$. 
\end{theorem}

\xproof{
Let $x: \mathbb{N}_a\rightarrow\mathbb{R}$ satisfy the initial conditions
	\[
	\left\{\begin{array}{l}x(a)=A, \\ x(a+1)=A+B.\end{array}\right.
	\]
Furthermore, for $t\in\N_{a+1}$, let $x(t+1)$ satisfy the recursive equation
	\begin{align*}
	x(t+1)&=x(t)-\sum_{s=a+1}^{t}\frac{(t-s+2)^{\overline{-\nu}}}{\Gamma(1-\nu)}\nabla x(s)\\
	&\quad + \frac{1}{p(t+1)}\left[h(t)-q(t)x(t)+p(t)\sum_{\tau=a+1}^{t}\frac{(t-\tau+1)^{\overline{-\nu}}}{\Gamma(1-\nu)}\nabla x(\tau)\right].
	\end{align*}
Note that as defined, $x(t+1)$ is uniquely determined from the values of $x(a)$, $x(a+1)$, \ldots, $f(t-1)$, $f(t)$, for $t \in \N_{a+1}$. Furthermore,
	\begin{align*}
	&\nabla x(t+1)+\sum_{s=a+1}^{t}\frac{(t-s+2)^{\overline{-\nu}}}{\Gamma(1-\nu)}\nabla x(s)\\
	&=\frac{1^{\overline{-\nu}}}{\Gamma(1-\nu)}\nabla x(t+1)+\sum_{s=a+1}^{t}\frac{(t-s+2)^{\overline{-\nu}}}{\Gamma(1-\nu)}\nabla x(s)\\
	&=\sum_{s=a+1}^{t+1}\frac{(t-s+2)^{\overline{-\nu}}}{\Gamma(1-\nu)}\nabla x(s)\\
	&=\nabla_a^{-(1-\nu)}\nabla x(t+1)\\
	&=\nabla_{a*}^{\nu}x(t+1).\\
	\end{align*}
So we have that 
\[
	\nabla_{a*}^{\nu}x(t+1) =\frac{1}{p(t+1)}\left[h(t)-q(t)x(t)+p(t)\sum_{\tau=a+1}^{t}\frac{(t-\tau+1)^{\overline{-\nu}}}{\Gamma(1-\nu)}\nabla x(\tau)\right].
\]
Then
	\begin{align*}
	&p(t+1)\nabla_{a*}^{\nu}x(t+1)\\
	&=h(t)-q(t)x(t)+p(t)\sum_{\tau=a+1}^{t}\frac{(t-\tau+1)^{\overline{-\nu}}}{\Gamma(1-\nu)}\nabla x(\tau)\\
	&=h(t)-q(t)x(t)+p(t)\nabla_a^{-(1-\nu)}\nabla x(t)\\
	&=h(t)-q(t)x(t)+p(t)\nabla_{a*}^{\nu}x(t),
	\end{align*}
which implies 
	\begin{equation*}
	\nabla[p(t+1)\nabla_{a*}^{\nu}x(t+1)]=h(t)-q(t)x(t).
	\end{equation*}
Thus, by rearranging, we have that
	\begin{equation*}
	\nabla[p(t+1)\nabla_{a*}^{\nu}x(t+1)]+q(t)x(t)=h(t).
	\end{equation*}
Therefore, for any value of $t\in\mathbb{N}_{a+1}$, $x(t)$ satisfies the IVP, So a solution exists. Reversing the preceding algebraic steps shows that if some function $y(t)$ is a solution to the IVP, it must be the same solution as our original $x(t)$. Therefore a unique solution exists.
}

The following lemma shows that initial conditions behave nicely when dealing with the Caputa nabla fractional difference.

\begin{lemma} \label{whocap}
Let $0<\nu<1$ and let $x:\mathbb{N}_{a}\rightarrow\mathbb{R}$. Then
	\[
	\nabla_{a*}^{\nu}x(a+1)=\nabla x(a+1)
	\]
for $t\in\mathbb{N}_{a+1}$. 
\end{lemma}

\xproof{
Let $0<\nu<1$.  Then by definition of the Caputo difference,  
	\begin{align*}
	\nabla_{a*}^{\nu}x(a+1)&=\nabla_a^{-(1-\nu)}\nabla x(a+1)\\
	&=\sum_{\tau=a+1}^{a+1} \frac{(a+1-\tau+1)^{\overline{-\nu}}}{\Gamma(1-\nu)}\nabla x(\tau)\\
	&=\nabla x(a+1),
	\end{align*}
for $t \in \N_{a+1}$.
} 
 
The next theorem and corollary show that the self-adjoint fractional nabla difference equation behaves very similar to a second order difference equation.

\begin{theorem} [General Solution of the Homogeneous Equation] \label{genhom}
Suppose $x_1,x_2:\mathbb{N}_a\rightarrow\mathbb{R}$ are linearly independent solutions to $L_ax(t)=0$. Then the general solution to $L_ax(t)=0$ is given by 
	\[
	x(t)=c_1x_1(t)+c_2x_2(t),
	\]
for $t\in\N_{a+1}$, where $c_1,c_2\in\mathbb{R}$ are arbitrary constants. 
\end{theorem}

\xproof{
Let $x_1,x_2:\mathbb{N}_a\rightarrow\mathbb{R}$ be linearly independent solutions of $L_ax(t)=0$. Then there exist constants $\alpha,\beta,\gamma,\delta\in\mathbb{R}$ for which $x_1,x_2$ are the unique solutions to the IVPs
	\[
	\begin{array}{lcr}
\left\{\begin{array}{l}Lx_1=0, \quad t \in \N_{a+1},\\x_1(a)=\alpha,\\ \nabla x_1(a+1)=\beta,\end{array}\right.
&\mathrm{and}&
\left\{\begin{array}{l}Lx_2=0, \quad t \in \N_{a+1},\\x_2(a)=\gamma,\\ \nabla x_2(a+1)=\delta.\end{array}\right.
	\end{array}
	\]
Since $L_a$ is a linear operator, we have for any $c_1,c_2\in\mathbb{R}$
	\[
	L_a[c_1x_1(t)+c_2x_2(t)]=c_1L_ax_1(t)+c_2L_ax_2(t)=0,
	\]
so $x(t)=c_1x_1(t)+c_2x_2(t)$ solves $L_ax(t)=0$. Conversely, suppose $x:\mathbb{N}_a\rightarrow\mathbb{R}$ solves $L_ax(t)=0$. Note that $x(t)$ solves the IVP
	\[
	\begin{cases}
	&Lx=0, \quad t \in \N_{a+1}, \\
	&x(a)=A, \\
	&\nabla x(a+1)=B,
	\end{cases}
	\]
for some $A,B \in \R$.
We show that the matrix equation
\begin{equation} \label{eq:mat}
\left[\begin{array}{cc}x_1(a) & x_2(a) \\ \nabla x_1(a+1) & \nabla x_2(a+1) \end{array}\right]\left[\begin{array}{c}c_1 \\ c_2\end{array}\right] = \left[\begin{array}{cc}A\\B\end{array}\right]
\end{equation}
has a unique solution for $c_1,c_2$. The above matrix equation can be equivalently expressed as 
\[
\left[\begin{array}{cc}\alpha & \gamma \\ \beta & \delta \end{array}\right]\left[\begin{array}{c}c_1 \\ c_2\end{array}\right] = \left[\begin{array}{cc}A\\B\end{array}\right].
\]
Suppose by way of contradiction that
\[
\left|\begin{array}{cc}\alpha & \gamma \\ \beta & \delta \end{array}\right|=0.
\]
Then, without loss of generality, there exists a constant $k\in\mathbb{R}$ for which $\alpha=k\gamma$ and $\beta=k\delta$. Then $x_1(a) = \alpha =k\gamma=kx_2(a)$, and $\nabla x_1(a+1)= \beta = k\delta=k \nabla x_2(a+1)$. Since $kx_2(t)$ solves $L_ax(t)=0$, we have that $x_1(t)$ and $kx_2(t)$ solve the same IVP.  By uniqueness, $x_1(t)=kx_2(t)$. But then $x_1(t)$ and $x_2(t)$ are linearly dependent, which is a contradiction. Therefore, the matrix equation \eqref{eq:mat} must have a unique solution, so $x(t)$ and $c_1x_1(t) + c_2x_2(t)$ solve the same IVP, and so by uniqueness in Theorem \ref{exuq}, every solution to $L_a x(t) = 0$ can be uniquely expressed as a linear combination of $x_1(t)$ and $x_2(t)$.
}

\begin{corollary}[General Solution of the Nonhomogeneous Equation]\label{gennon}
Suppose $x_1,x_2:\mathbb{N}_a\rightarrow\mathbb{R}$ are linearly independent solutions of $L_ax(t)=0$ and $y_0:\mathbb{N}_a\rightarrow\mathbb{R}$ is a particular solution to $L_ax(t)=h(t)$ for some $h:\mathbb{N}_{a+1}\rightarrow\mathbb{R}$. Then the general solution of $L_ax(t)=h(t)$ is given by 
	\[
	x(t)=c_1x_1(t)+c_2x_2(t)+y_0(t),
	\]
for $t\in\N_{a+1}$ and where $c_1,c_2\in\mathbb{R}$ are arbitrary constants. 
\end{corollary}

\xproof{
Since $L_a$ is a linear operator, one can show that $x(t)=c_1x_1(t)+c_2x_2(t)+y_0(t)$ solves $L_ax(t)=h(t)$ for any $c_1,c_2\in\mathbb{R}$ in a similar way as in Theorem \ref{genhom}. Conversely, suppose $x:\mathbb{N}_a\rightarrow\mathbb{R}$ solves $L_ax(t)=h(t)$. Again note that $x(t)$ solves the IVP
	\[
	\begin{cases}
	&Lx=h(t),\quad t \in \N_{a+1}, \\
	&x(a)=A, \\
	&x(a+1)=B,
	\end{cases}
	\]
for some $A,B \in \R$. Since $y_0(t)$ is a particular solution of $L_ax(t)=h(t)$, there exist unique constants $c_1,c_2\in\mathbb{R}$ for which $x_h(t)=c_1x_1(t)+c_2x_2(t)$ solves the IVP
	\[
	\begin{cases}
	&Lx_h=0,\\
	&x_h(a)=A-y_0(a),\\
	&x_h(a+1)=B-y_0(a+1),
	\end{cases}
	\]
for $t\in\mathbb{N}_{a+1}$. Observe that $x_h(t)+y_0(t)$ satisfies $L_ax(t)=h(t)$. Further, 
	\[
	x_h(a)+y_0(a)=A-y_0(a)+y_0(a)=A,
	\]
and 
	\[
	x_h(a+1)+y_0(a+1)=B-y_0(a+1)+y_0(a+1)=B.
	\]
Therefore $x(t)$ and $x_h(t)+y_0(t)$ solve the same IVP.  Then by uniqueness of IVPs, $x(t)=c_1x_1(t)+c_2x_2(t)+y_0(t)$, and thus any solution to $L_ax(t)=h(t)$ may be written in this form. 
}

\section{Initial value problems for the fractional self-adjoint equation}
In this section we develop techniques to solve initial value problems for the fractional self-adjoint operator involving the Caputo difference. See Brackins \cite{Brackins} for a similar development using the Riemann-Liouville definition of a fractional difference.

\begin{definition}\label{caudef}
The \textit{Cauchy function for $L_ax(t)$} is the function $x:\mathbb{N}_s\times\mathbb{N}_a\rightarrow\mathbb{R}$ that satisfies the IVP
\begin{equation} \label{eq:cauwho}
\left\{\begin{array}{l}L_sx(t,s)=0, \quad t \in \N_{s+1},\\x(s,s)=0,\\ \nabla x(s+1,s)=\frac{1}{p(s+1)},\end{array}\right.
\end{equation}
for any fixed $s \in \N_a$.
\end{definition}

\begin{remark}\label{altcaudef}
Note by Lemma \ref{whocap}, the IVP \eqref{eq:cauwho} is equivalent to the IVP
\[
\left\{\begin{array}{l}L_sx(t,s)=0, \quad t \in \N_{s+1} \\x(s,s)=0,\\ \nabla_{s*}^{\nu} x(s+1,s)=\frac{1}{p(s+1)}
\end{array}\right.
\]
for any fixed $s \in \N_a$.
\end{remark}

\begin{theorem}[Variation of Constants] \label{voc}
Let $h: \N_{a+1} \to \R$. Then solution to the IVP
	\[
	\begin{cases}
	&L_ay(t)=h(t), \quad t \in \N_{a+1},\\
	&y(a) =0,\\
	&\nabla y(a+1)=0,
	\end{cases}
	\]
is given by 
	\[y(t)=\int_a^t x(t,s)h(s)\nabla s,
	\]
where $x(t,s)$ is the Cauchy function for the homogenous equation and where $y : \N_a \to \R$.
\end{theorem}

\xproof{
Note that
\[
y(a)=\sum_{s=a+1}^a x(a,s)h(s)=0,
\]
so the first initial condition holds.
By the definition of the Caputo difference, 
	\[
	\nabla_{s*}^{\nu}x(t+1,s)=\nabla_s^{-(1-\nu)}\nabla x(t+1,s)=\sum_{\tau=s+1}^{t+1} \frac{(t+1-\tau+1)^{\overline{-\nu}}}{\Gamma(-\nu)}\nabla x(\tau,s).
	\]
Then
	\begin{align*}
	&\nabla[p(t+1)\nabla_{s*}^{\nu}x(t+1,s)] \\
	&= p(t+1)\sum_{\tau=s+1}^{t+1} \frac{(t-\tau+2)^{\overline{-\nu}}}{\Gamma(1-\nu)}\nabla x(\tau,s)-p(t)\sum_{\tau=s+1}^{t} \frac{(t-\tau+1)^{\overline{-\nu}}}{\Gamma(1-\nu)}\nabla x(\tau,s)\\
	&=p(t+1)\sum_{\tau=s+1}^{t+1} \frac{(t-\tau+2)^{\overline{-\nu}}}{\Gamma(1-\nu)}\nabla x(\tau,s)-p(t)\sum_{\tau=s+1}^{t+1} \frac{(t-\tau+1)^{\overline{-\nu}}}{\Gamma(1-\nu)}\nabla x(\tau,s)\\
	&\quad +p(t)\frac{(t-t-1+1)^{\overline{-\nu}}}{\Gamma(1-\nu)}\nabla x(t+1,s)\\
	&=\sum_{\tau=s+1}^{t+1}\frac{(t-\tau+2)^{\overline{-\nu}}p(t+1)-(t-\tau+1)^{\overline{-\nu}}p(t)}{\Gamma(1-\nu)}\nabla x(\tau,s).
	\end{align*}
So we have that
\begin{multline}\label{varconintermediate}
\nabla[p(t+1)\nabla_{s*}^{\nu}x(t+1,s)] \\
\quad= \sum_{\tau=s+1}^{t+1}\frac{(t-\tau+2)^{\overline{-\nu}}p(t+1)-(t-\tau+1)^{\overline{-\nu}}p(t)}{\Gamma(1-\nu)}\nabla x(\tau,s).
\end{multline}
Now we consider $y(t)$.  Note that
	\begin{align*}
	\nabla y(t) &= \sum_{s=a+1}^{t}x(t,s)h(s) -\sum_{s=a+1}^{t-1}x(t-1,s)h(s)\\
	&= \sum_{s=a+1}^{t-1}x(t,s)h(s) +x(t,t)h(t)-\sum_{s=a+1}^{t-1}x(t-1,s)h(s)\\
	&= \sum_{s=a+1}^{t-1}\nabla x(t,s)h(s).
	\end{align*}
From here, observe that
\[
\nabla y(a+1) = \sum_{s=a+1}^{a+1-1}x(t,s)h(s)=0,
\]
so the second initial condition holds.
Then by the definition of the Caputo difference, 
	\begin{align*}
	\nabla_{a*}^{\nu}y(t)&=\nabla_a^{-(1-\nu)}\nabla y(t)\\
	&=\sum_{\tau=a+1}^{t} \frac{(t-\tau+1)^{\overline{-\nu}}}{\Gamma(1-\nu)}\nabla y(\tau)\\
	&=\sum_{\tau=a+1}^{t} \frac{(t-\tau+1)^{\overline{-\nu}}}{\Gamma(1-\nu)} \sum_{s=a+1}^{\tau-1}\nabla x(\tau,s)h(s)\\
	&=\sum_{\tau=a+1}^{t} \sum_{s=a+1}^{\tau-1} \frac{(t-\tau+1)^{\overline{-\nu}}}{\Gamma(1-\nu)}\nabla x(\tau,s)h(s).
	\end{align*}
Then from (\ref{varconintermediate})
	\begin{align*}
	&\nabla[p(t+1)\nabla_{a*}^{\nu}y(t+1)] \\
	&= p(t+1)\sum_{\tau=a+1}^{t+1} \sum_{s=a+1}^{\tau-1} \frac{(t-\tau+2)^{\overline{-\nu}}}{\Gamma(1-\nu)}\nabla x(\tau,s)h(s) \\
	&\quad - p(t)\sum_{\tau=a+1}^{t} \sum_{s=a+1}^{\tau-1} \frac{(t-\tau+1)^{\overline{-\nu}}}{\Gamma(1-\nu)}\nabla x(\tau,s)h(s)\\
	&=\sum_{\tau=a+1}^{t+1}\sum_{s=a+1}^{\tau-1} \frac{(t-\tau+2)^{\overline{-\nu}}p(t+1)-(t-\tau+1)^{\overline{-\nu}}p(t)}{\Gamma(1-\nu)}\nabla x(\tau,s)h(s)\\
	&=\sum_{s=a+1}^{t}\sum_{\tau=s+1}^{t+1} \frac{(t-\tau+2)^{\overline{-\nu}}p(t+1)-(t-\tau+1)^{\overline{-\nu}}p(t)}{\Gamma(1-\nu)}\nabla x(\tau,s)h(s)\\
	&=p(t+1)\nabla x(t+1,t)h(t)\\
	&\quad + \sum_{s=a+1}^{t-1}\sum_{\tau=s+1}^{t+1} \frac{(t-\tau+2)^{\overline{-\nu}}p(t+1)-(t-\tau+1)^{\overline{-\nu}}p(t)}{\Gamma(1-\nu)}\nabla x(\tau,s)h(s)\\
	&=p(t+1)\frac{h(t)}{p(t+1)}+\sum_{s=a+1}^{t-1}\nabla[p(t+1)\nabla_{s*}^{\nu}x(t+1,s)]h(s)\\
	&=h(t)+\sum_{s=a+1}^{t-1}\nabla[p(t+1)\nabla_{s*}^{\nu}x(t+1,s)]h(s).
	\end{align*}
Therefore
	\begin{align*}
	L_ay(t)&=\nabla[p(t+1)\nabla_{a*}^{\nu}y(t+1)]+q(t)y(t)\\
	&=h(t)+\sum_{s=a+1}^{t-1}\nabla[p(t+1)\nabla_{s*}^{\nu}x(t+1,s)]h(s)\\
	&\quad+\sum_{s=a+1}^{t-1}q(t)x(t,s)h(s)+q(t)x(t,t)h(t)\\
	&=h(t)+\sum_{s=a+1}^{t-1}\bigg[\nabla[p(t+1)\nabla_{s*}^{\nu}x(t+1,s)]+q(t)x(t,s)\bigg]h(s)\\
	&=h(t)+\sum_{s=a+1}^{t-1} L_sx(t,s)h(s)\\
	&=h(t).
	\end{align*}
Thus $y:\N_a \to \R$ solves the IVP for $t \in \N_{a+1}$. 
}

\begin{theorem}[Variation of Constants with Non-Zero Initial Conditions]\label{vocnon}
The solution to the IVP
	\[
	\begin{cases}
	&L_a y(t)=h(t), \quad t \in \N_{a+1}, \\
	&y(a)=A,\\
	&\nabla y(a+1)=B,
	\end{cases}
	\]
where $A,B\in \R$ are arbitrary constants, is given by 
	\[
	y(t)=y_0(t)+\int_a^t x(t,s)h(s)\nabla s,
	\]
where $y_0(t)$ solves the IVP
	\[
	\begin{cases}
	&L_a y_0(t)=0, \quad t \in \N_{a+1},\\
	&y_0(a)=A,\\
	&\nabla y_0(a+1)=B.
	\end{cases}
	\]	
\end{theorem}

\xproof{
The proof follows from Theorem \ref{voc} by linearity.
}

\begin{example} \label{exa}
Find the Cauchy function for 
	\[
	\nabla[p(t+1)\nabla_{a*}^{\nu}y(t+1)]=0, \quad t\in \N_{a+1}.
	\]
Consider $\nabla[p(t+1)\nabla_{s*}^{\nu}x(t+1,s)]=0$.  Integrating both sides from $s$ to $t$ and applying the Fundamental Theorem of Nabla Calculus along with the second initial condition from Remark \ref{altcaudef} yields
	\begin{align*}
	p(t+1)\nabla_{s*}^{\nu}x(t+1,s) - p(s+1)\nabla_{s*}^{\nu}x(s+1,s)&=0\\
	p(t+1)\nabla_{s*}^{\nu}x(t+1,s) - 1&=0 \\
	\nabla_{s*}^\nu x(t+1,s) &= \frac{1}{p(t+1)}.
	\end{align*}
By the definition of the Caputo difference, this is equivalent to
	\begin{align*}
	\nabla_{s}^{-(1-\nu)} \nabla x(t+1,s) &=\frac{1}{p(t+1)}\\
	\nabla_s^{1-\nu}\nabla_{s}^{-(1-\nu)}\nabla x(t+1,s)&=\nabla_s^{1-\nu} \frac{1}{p(t+1)}\\
	\nabla x(t+1,s) &=\nabla_s^{1-\nu} \frac{1}{p(t+1)},
	\end{align*}
after applying a composition rule from Theorem \ref{comprules}.
Replacing $t+1$ with $t$ yields
	\[
	\nabla x(t,s)=\nabla_s^{1-\nu} \frac{1}{p(t)}.
	\]
Applying the Fundamental Theorem after integrating both sides from $s$ to $t$ and applying the first initial condition from Remark \ref{altcaudef} yields
	\begin{align*}
	x(t,s)-x(s,s) &=\int_s^t \nabla_s^{1-\nu} \frac{1}{p(\tau)}\nabla \tau\\
	x(t,s)&=\int_s^t \nabla\nabla_s^{-\nu} \frac{1}{p(\tau)}\nabla \tau\\
	&=\bigg[\nabla_s^{-\nu} \frac{1}{p(\tau)}\bigg]_{\tau=s}^{\tau=t}\\
	&=\nabla_s^{-\nu} \frac{1}{p(t)}-\nabla_s^{-\nu} \frac{1}{p(s)}\\
	&=\nabla_s^{-\nu} \frac{1}{p(t)}.
	\end{align*}
Therefore the Cauchy function is 
	\[
	x(t,s)=\nabla_s^{-\nu}\frac{1}{p(t)}=\sum_{\tau=s+1}^t \frac{(t-\tau+1)^{\overline{\nu-1}}}{\Gamma(\nu)} \left(\frac{1}{p(\tau)}\right).
	\]
\end{example}

\begin{example} \label{exb}
Find the Cauchy function for 
	\[
	\nabla\nabla_{a*}^{\nu}y(t+1)=0, \quad t \in \N_{a+1}.
	\]
Notice that this is a particular case of the previous example, where $p(t)\equiv 1$.  Then the Cauchy function is 
	\begin{align*}
	x(t,s)&=(\nabla_s^{-\nu} 1)(t)\\
	&=\int_s^t \frac{(t-\tau+1)^{\overline{\nu-1}}}{\Gamma(\nu)}\nabla \tau\\
	&=-\frac{(t-\tau)^{\overline{\nu}}}{\Gamma(\nu+1)}\bigg|_{\tau=s}^{\tau=t}\\
	&=\frac{(t-s)^{\overline{\nu}}}{\Gamma(\nu+1)}.
	\end{align*}
Also note that if you take $\nu=1$ as in the whole order self-adjoint case, the Cauchy function simplifies to $x(t,s)=t-s$.
\end{example}

\begin{example} \label{exc}
Find the solution to the IVP
	\[
	\begin{cases}
	&\nabla[p(t+1)\nabla_{a*}^{\nu}y(t+1)]=h(t), \quad t \in \N_{a+1}, \\
	&y(a)= 0, \\
	&\nabla y(a+1)=0.
	\end{cases}
	\]
From Theorem \ref{voc}, we know that the solution is given by
	\[
	y(t)=\int_a^t x(t,s)h(s)\nabla s.
	\]
By Example \ref{exa}, we know the Cauchy function for the above difference equation is $x(t,s)=\nabla_s^{-\nu}\frac{1}{p(t)}$. 
Then the solution is given by
	\[
	y(t)=\int_a^t \nabla_s^{-\nu} \frac{1}{p(t)} h(s) \nabla s.
	\]
\end{example}

\begin{example}
Solve the IVP
	\[
	\begin{cases}
	&\nabla \nabla_{0*}^{0.6} x(t+1)=t, \quad t\in \N_1,\\
	&x(0)= 0, \\
	&\nabla x(1)=0.
	\end{cases}
	\]
	This is a particular case of Example \ref{exc} where $h(t)=t$, $a=0$, and $\nu=0.6$. Then
	
	\begin{align*}
	y(t)&=\int_0^t x(t,s)s\nabla s\\
	&=\sum_{s=1}^t \frac{1}{\Gamma(1.6)}(t-s)^{\overline{0.6}}s,\\
	\end{align*}
and after summing by parts and applying Theorem \ref{fund}, we get that the solution is
	\[
	y(t)=\frac{(t-1)^{\overline{2.6}}}{\Gamma(3.6)}
	\]
\end{example}

\section{Boundary value problems of the fractional self-adjoint equation}
In this section we develop techniques to solve boundary value problems for the fractional self-adjoint operator involving the Caputo difference. See Brackins \cite{Brackins} for a similar development using the Riemann-Liouville definition of a fractional difference.

We are interested in the boundary value problems (BVPs)
	\begin{equation}	\label{trivbvp}
	\begin{cases}
	&L_a x(t)=0, \quad t \in \N_{a+1}^{b-1}, \\ 
	&\alpha x(a)-\beta \nabla x(a+1)=0, \\ 
	&\gamma x(b) + \delta\nabla x(b)=0,
	\end{cases}
	\end{equation}
and
	\begin{equation} \label{nontrivbvp}
	\begin{cases}
	&L_a x(t)=h(t), \quad t \in \N_{a+1}^{b-1},\\ 
	&\alpha x(a)-\beta \nabla x(a+1)=A, \\ 
	&\gamma x(b) + \delta\nabla x(b)=B,
	\end{cases}
	\end{equation}
where $h:\mathbb{N}_{a+1}^{b-1}\rightarrow\mathbb{R}$ and $\alpha,\beta,\gamma,\delta,A,B\in\mathbb{R}$ for which $\alpha^2+\beta^2>0$ and $\gamma^2+\delta^2>0$. Note that despite the fact that the difference equations above hold for $t \in \N_{a+1}^{b-1}$, the solution $x(t)$ for each BVP is defined on the domain of $\N_a^b$. We are primarily interested in cases where the BVP $(\ref{trivbvp})$ has only the trivial solution. 

\begin{theorem} \label{bvpuniq}
Assume $(\ref{trivbvp})$ has only the trivial solution. Then $(\ref{nontrivbvp})$ has a unique solution. 
\end{theorem}

\xproof{
Let $x_1,x_2:\mathbb{N}_a\rightarrow\mathbb{R}$ be linearly independent solutions to $L_a x(t)=0$. By Theorem 13, a general solution to $L_a x(t)=0$ is given by 
	\[
	x(t)=c_1x_1(t)+c_2x_2(t),
	\]
where $c_1,c_2\in\mathbb{R}$ are arbitrary constants. If $x(t)$ solves the boundary conditions in $(\ref{trivbvp})$, then $x(t)$ is the trivial solution, which is true if and only if $c_1=c_2=0$.  This is true if and only if the system of equations 
\[
\left\{\begin{array}{l}\alpha[c_1x_1(a)+c_2x_2(a)]-\beta\nabla_{a*}^{\nu}[c_1x_1(a+1)+c_2x_2(a+1)]=0,\\
\gamma[c_1x_1(b)+c_2x_2(b)]+\delta\nabla_{a*}^{\nu}[c_1x_1(b)+c_2x_2(b)]=0,\end{array}\right.
\]
or equivalently, 
\[
\left\{\begin{array}{l}c_1[\alpha x_1(a)-\beta\nabla_{a*}^{\nu}x_1(a+1)]+c_2[\alpha x_2(a)-\beta\nabla_{a*}^{\nu}x_2(a+1)]=0,
\\ c_1[\gamma x_1(b)+\delta\nabla_{a*}^{\nu}x_1(b)]+c_2[\gamma x_2(b)+\delta\nabla_{a*}^{\nu}x_2(b)]=0,\end{array}\right.
\]
has only the trivial solution  where $c_1=c_2=0$. In other words, $x(t)$ solves $(\ref{trivbvp})$ if and only if 
\begin{align*}
D&:=\left|\begin{array}{cc}\alpha x_1(a)-\beta\nabla_{a*}^{\nu}x_1(a+1) & \alpha x_2(a)-\beta\nabla_{a*}^{\nu}x_2(a+1) \\
\gamma x_1(b) + \delta\nabla_{a*}^{\nu}x_1(b) & \gamma x_2(b)+\delta\nabla_{a*}^{\nu}x_2(b) \end{array}\right| \\
&\neq 0.
\end{align*}
Now consider $(\ref{nontrivbvp})$. By Corollary \ref{gennon}, a general solution to $L_a y(t)=h(t)$ is 
	\[
	y(t)=a_1x_1(t)+a_2x_2(t)+y_0(t),
	\]
where $a_1,a_2\in\mathbb{R}$ are arbitrary constants and $y_0:\mathbb{N}_a\rightarrow\mathbb{R}$ is a particular solution of $L_a y(t)=h(t)$. Consider the system of equations 
\[
\left\{\begin{array}{l}a_1[\alpha x_1(a)-\beta\nabla_{a*}^{\nu}x_1(a+1)]+a_2[\alpha x_2(a)-\beta\nabla_{a*}^{\nu}x_2(a+1)]\\
\quad =A-\alpha y_0(a)+\beta\nabla_{a*}^{\nu}y_0(a+1),\\
a_1[\gamma x_1(b) + \delta\nabla_{a*}^{\nu}x_1(b)] + a_2[\gamma x_2(b)+\delta\nabla_{a*}^{\nu}x_2(b)]=B-\gamma y_0(b)-\delta\nabla_{a*}^{\nu}y_0(b),\end{array}\right.
\]
for arbitrary $A,B\in\mathbb{R}$ as in $(\ref{nontrivbvp})$. Since $D\neq0$, this system has a unique solution for $a_1,a_2$. It may be shown algebraically that this system is equivalent to
\[
\left\{\begin{array}{l}\alpha[a_1x_1(a)+a_2x_2(a)+y_0(a)]-\beta\nabla_{a*}^{\nu}[a_1x_1(a+1)+a_2x_2(a+1)+y_0(a+1)]\\\quad=A,\\
\gamma[a_1x_1(b)+a_2x_2(b)+y_0(b)]+\delta\nabla_{a*}^{\nu}[a_1x_1(b)+a_2x_2(b)+y_0(b)]=B,\end{array}\right.
\]
so $y(t)$ satisfies the boundary conditions for $(\ref{nontrivbvp})$. Therefore for any $A,B\in\mathbb{R}$, $(\ref{nontrivbvp})$ has a unique solution.
}

\begin{theorem}
Let $$\rho := \alpha \gamma \nabla^{-\nu}_a \frac{1}{p(b)} + \frac{\alpha\delta}{p(b)} + \frac{\beta\gamma}{p(a+1)}.$$ Then the BVP
	\[
	\begin{cases}
	&\nabla[p(t+1) \capnab x(t+1)] = 0, \quad t\in\mathbb{N}_{a+1}^{b-1},\\
	&\alpha x(a) - \beta\capnab x(a+1) = 0, \\
	&\gamma x(b) + \delta\capnab x(b) = 0,
	\end{cases}
	\]
has only the trivial solution if and only if $\rho \neq 0$.
\end{theorem}

\xproof{
Note that $x_1(t) = 1, x_2(t) = \nabla^{-\nu}_a \frac{1}{p(t)}$ are linearly independent solutions to 
	\[
	\nabla [p(t+1) \capnab x(t+1)] = 0.
	\]
Then a general solution of the difference equation is given by
	\[
	x(t) = c_1 x_1(t) + c_2x_2(t) = c_1 + c_2 \nabla^{-\nu}_a \frac{1}{p(t)}.
	\]
Consider the boundary conditions $\alpha x(a) - \beta\capnab x(a+1) = 0,$ and $\gamma x(b) + \delta\capnab x(b) = 0$. These boundaries give us
	\[
	c_1\alpha + c_2\bigg[- \frac{\beta}{p(a+1)}\bigg] = 0,
	\] 
	\[
	c_1\gamma+ c_2\bigg(\frac{\delta}{p(b)} + \gamma \nabla^{-\nu}_a \frac{1}{p(b)}\bigg) = 0.
	\] 
Converting this into a linear system yields
\[\left[\begin{array}{cc}
\alpha & - \frac{\beta}{p(a+1)} \\ \gamma & \frac{\delta}{p(b)} + \gamma \nabla^{-\nu}_a \frac{1}{p(b)} \end{array}\right] 
\left[\begin{array}{c}c_1\\c_2\end{array}\right]=\left[\begin{array}{c}0\\0\end{array}\right].\]
Consider the determinant of the coefficient matrix, 
\[\left|\begin{array}{cc}
\alpha & - \frac{\beta}{\rho(a+1)} \\ \gamma & \frac{\delta}{p(b)} + \gamma \nabla^{-\nu}_a \frac{1}{p(b)} \end{array}\right| = \alpha \gamma \nabla^{-\nu}_a \frac{1}{p(b)} + \frac{\alpha\delta}{p(b)} + \frac{\beta\gamma}{p(a+1)} = \rho.
\]
By properties of invertible matrices, the BVP has only the trivial solution if and only if $\rho \neq 0$.
}

\begin{definition}
Assume that (\ref{trivbvp}) has only the trivial solution. Then we define the \textit{Green's function for the homogeneous BVP (\ref{trivbvp}}), $G(t,s)$, by
\[
G(t,s):=\left\{\begin{array}{lcr}u(t,s),&\quad&a\le t\le s\le b, \\
v(t,s),&\quad&a\le s\le t\le b,\end{array}\right.
\]
where $u(t,s)$ solves the BVP
\[
\left\{\begin{array}{l}L_a u(t)=0, \quad t \in \N_{a+1}^{b-1},
\\ \alpha u(a,s)-\beta\nabla u(a+1,s)=0,\\ \gamma u(b,s)+\delta\nabla u(b,s)=-[\gamma x(b,s)+\delta \nabla x(b,s)],\end{array}\right.
\]
for each fixed $s\in\mathbb{N}_a^b$ and where $x(t,s)$ is the Cauchy function for $L_a x(t)$. Then we define
	\[
	v(t,s):=u(t,s)+x(t,s).
	\]
\end{definition}

\begin{theorem}[Green's Function Theorem] \label{grnfxnthm}
If $(\ref{trivbvp})$ has only the trivial solution, then the solution to $(\ref{nontrivbvp})$ where $A=B=0$ is given by 
	\[
	y(t)=\int_a^bG(t,s)h(s)\nabla s,
	\]
where the $G(t,s)$ is the Green's function for the homogeneous BVP (\ref{trivbvp}). \end{theorem}

\xproof{
First note that by Theorem \ref{bvpuniq}, $u(t,s)$ for each fixed $s\in\mathbb{N}_{a}^{b}$ is well-defined.
Let 
	\begin{align*}
	y(t) = \int^b_a G(t,s) h(s) \nabla s &= \int^t_a G(t,s) h(s) \nabla s + \int^b_t G(t,s) h(s) \nabla s\\
    &= \int^t_a v(t,s) h(s) \nabla s + \int^b_t u(t,s) h(s) \nabla s\\
    &= \int^t_a [u(t,s) + x(t,s)] h(s) \nabla s + \int^b_t u(t,s) h(s) \nabla s\\
    &= \int^b_a u(t,s) h(s) \nabla s + \int^t_a x(t,s) h(s) \nabla s\\
    &= \int^b_a u(t,s) h(s) \nabla s + z(t),\\
    \end{align*}
where $z(t):=\int^t_a x(t,s) h(s) \nabla s$.  Since $x(t,s)$ is the Cauchy function for $L_a x(t)=0$, by Theorem \ref{voc}, $z(t)$ solves the IVP
	\[
	\left\{
		\begin{array}{lr}
			L_a z(t) = h(t), \quad t \in \N_{a+1}^{b-1}, \\
			z(a) = 0,\\
			\nabla z(a+1) = 0.
		\end{array}
	\right.
	\]
Then,
\begin{align*}
L_a y(t) &= \int^b_a L_a u(t,s) h(s) \nabla s + L_a z(t)\\
           &= 0 + h(t) = h(t), 
\end{align*}
for $t\in\mathbb{N}_{a+1}^{b-1}$, so the difference equation is satisfied.  Now we check the boundary conditions.  At $t=a$, we have
	\[
	\alpha y(a) - \beta  \nabla y(a+1) = \int^b_a [\alpha u(a,s) - \beta\nabla u(a+1,s)] h(s) \nabla s + [\alpha z(a) - \beta \nabla z(a+1)] = 0,
	\]
and at $t=b$, we have
\begin{align*}
&\gamma y(b)+\delta \nabla y(b)\\
&=\gamma z(b)+\int_{a}^{b}\gamma u(b,s)h(s)\nabla s+\delta\nabla z(b)+\int_{a}^{b}\delta\nabla  u(b,s)h(s)\nabla s\\
&=\gamma \int_{a}^{b}x(b,s)h(s)\nabla s + \delta\nabla\int_{a}^{b} x(b,s)h(s)\nabla s +\int_{a}^{b}[\gamma u(b,s)+\delta\nabla u(b,s)]h(s)\nabla s\\
&=-\int_{a}^{b}[\gamma x(b,s)+\delta\nabla x(b,s)]h(s)\nabla s +\int_{a}^{b}[\gamma x(b,s)+\delta \nabla x(b,s)]h(s)\nabla s\\
&=0.
\end{align*}
}

\begin{corollary}
If $(\ref{trivbvp})$ has only the trivial solution, then the solution to $(\ref{nontrivbvp})$ with $A,B \in \R$ is given by 
	\[
	y(t)=z(t)+\int_a^bG(t,s)h(s)\nabla s,
	\]
where $z:\mathbb{N}_a^b\rightarrow\mathbb{R}$ is the unique solution to 
\[
\left\{\begin{array}{l}L_a z(t)=0,\\ \alpha z(a)-\beta\nabla z(a+1)=A,\\ \gamma z(b)+\delta\nabla z(b)=B.
\end{array}\right.
\]
\end{corollary}

\xproof{
This corollary follows directly from Theorem \ref{grnfxnthm} by linearity.
}

\begin{example}\label{greensbvp}
Find the Green's function for the boundary value problem
\[
\begin{cases}
\nabla[\nabla_{a*}^\nu y(t+1)] = 0, \quad t \in \N_{a+1}^{b-1},\\
y(a) = 0,\\
y(b) = 0.
\end{cases}
\]
The Green's function is given by 
\[
G(t,s) =
\begin{cases}
u(t,s), &a \le t \le s \le b,\\
v(t,s), &a \le s \le t \le b,
\end{cases}
\]
where $u(t,s)$, for each fixed $s\in \N_a^b$, solves the BVP
\[
\begin{cases}
\nabla[\nabla_{a*}^\nu u(t+1,s)] = 0,\quad t \in \N_{a+1}^{b-1}, \\
u(a,s) = 0,\\
u(b,s) = -x(b,s),
\end{cases}
\]
and $v(t,s) = u(t,s) + x(t,s)$.  By inspection, we find that $x_1(t)=1$ is a solution of 
\[
\nabla[\nabla_{a*}^\nu y(t+1)] = 0,
\]
for $t\in\N_{a+1}$.  Let $x_2(t)=(\nabla_a^{-\nu} 1)(t)$. Consider
\begin{align*}
\nabla[\nabla_{a*}^{\nu}x_2(t+1)] &= \nabla[\nabla_{a*}^{\nu}\nabla_a^{-\nu}(1)]\\
&= \nabla[\nabla_a^{-(N-\nu)}\nabla^N\nabla_a^{-\nu}(1)]\\
&= \nabla[\nabla_a^{-(N-\nu)}\nabla_a^{(N-\nu)}(1)]\\
&= \nabla [1]\\
&= 0,
\end{align*}
using Theorem \ref{comprules}. So we have that $x_2(t)$ solves $\nabla[\nabla_{a*}^\nu y(t+1)] = 0$.  Since $x_1(t)$ and $x_2(t)$ are linearly independent, by Theorem \ref{genhom}, the general solution is given by
\begin{align*}
y(t) = c_1 + c_2 (\nabla_a^{-\nu} 1)(t)
= c_1 + c_2 \frac{(t-a)^{\overline \nu}}{\Gamma(1+\nu)},
\end{align*}
and it follows that
\[
u(t,s) = c_1(s) + c_2(s)\frac{(t-a)^{\overline \nu}}{\Gamma(1+\nu)}.
\]
The boundary condition $u(a,s) = 0$ implies that $c_1(s) = 0$.  The boundary condition $u(b,s) = -x(b,s)$ then yields
\[
-x(b,s) = u(b,s) = c_2(s)\frac{(b-a)^{\overline \nu}}{\Gamma(1+\nu)}.
\]
From Example \ref{exb}, we know that
\[
x(b,s) = (\nabla_s^{-\nu} 1)(b) = \frac{(b-s)^{\overline \nu}}{\Gamma(1+\nu)} ,
\]
and thus
\[
c_2(s) = -\frac{(b-s)^{\overline \nu}}{(b-a)^{\overline \nu}}.
\]
Hence the Green's function is given by
\[
G(t,s) = \begin{cases}
-\dfrac{(b-s)^{\overline \nu}(t-a)^{\overline \nu}}{\Gamma(1+\nu)(b-a)^{\overline \nu}}, & a \le t \le s \le b,\\
-\dfrac{(b-s)^{\overline \nu}(t-a)^{\overline \nu}}{\Gamma(1+\nu)(b-a)^{\overline \nu}} + \dfrac{(t-s)^{\overline \nu}}{\Gamma(1+\nu)}, & a\le s \le t \le b.
\end{cases}
\]
\end{example}

\begin{remark}
Note that in the continuous and whole-order discrete cases, the Green's function is symmetric for the equivalent BVP in Example \ref{greensbvp}.  This is not necessarily true in the fractional case.  By way of counterexample, take $a = 0$, $b = 5$, and $\nu = 0.5$. Then computing we find that
\[
G(2,3) = u(2,3) = -\frac{(2)^{\overline{0.5}}(2)^{\overline{0.5}}}{\Gamma(1.5)(5)^{\overline{0.5}}} = -\frac{32}{35},
\]
but
\[
G(3,2) = v(3,2) = -\frac{(3)^{\overline{0.5}}(3)^{\overline{0.5}}}{\Gamma(1.5)(5)^{\overline{0.5}}} + \frac{(1)^{\overline{0.5}}}{\Gamma(1.5)} = -\frac{3}{7}.
\]
Thus for this particular BVP, unlike in the continuous and whole-order discrete cases, is not symmetric.
\end{remark}

\begin{theorem}\label{ineq}
The Green's function for the BVP
	\[
	\begin{cases}
	&\nabla\nabla_{a*}^{\nu}x(t+1)=0,\\
	&x(a)=x(b)=0,
	\end{cases}
	\]
for $t\in\N_{a+1}^{b-1}$, given by
	\[
	G(t,s)=
	\begin{cases}
	-\dfrac{(b-s)^{\overline \nu}(t-a)^{\overline \nu}}{\Gamma(1+\nu)(b-a)^{\overline \nu}}, & a \le t \le s \le b,\\
-\dfrac{(b-s)^{\overline \nu}(t-a)^{\overline \nu}}{\Gamma(1+\nu)(b-a)^{\overline \nu}} + \dfrac{(t-s)^{\overline \nu}}{\Gamma(1+\nu)}, & a\le s \le t \le b,
	\end{cases}
	\]
satisfies the inequalities
	\begin{enumerate}
	\item $\displaystyle G(t,s)\leq 0, $
	\item $\displaystyle G(t,s)\geq -\bigg(\frac{b-a}{4}\bigg)\bigg(\frac{\Gamma(b-a+1)}{\Gamma(\nu+1)\Gamma(b-a+\nu)}\bigg), $
	\item $\displaystyle \int_a^b |G(t,s)|\nabla s \leq \frac{(b-a)^2}{4\Gamma(\nu+2)}, $ 
	\newcounter{enumTemp}
    \setcounter{enumTemp}{\theenumi}
	\end{enumerate}
for $t\in\N_a^b$, and
	\begin{enumerate}
	 \setcounter{enumi}{\theenumTemp}
	\item $\displaystyle \int_a^b |\nabla G(t,s)|\nabla s \le \frac{b-a}{\nu+1},$
	\end{enumerate}
for $t\in\N_{a+1}^b$.
\end{theorem}

\xproof{
(1) Let $a\leq t\leq s\leq b$.  Then
	\[
	G(t,s)=u(t,s)=-\frac{(t-a)^{\overline{\nu}}(b-s)^{\overline{\nu}}}{\Gamma(\nu+1)(b-a)^{\overline{\nu}}}\leq 0,
	\]
	for each fixed $s \in \N_a^b$.
Now let $a\leq s<t\leq b$.  Then $G(t,s)=v(t,s)$, so we wish to show that $v(t,s)$ is non-positive.  First, we show that $v(t,s)$ is increasing. Taking the nabla difference with respect to $t$ yields
	\[
	\nabla_t \bigg[ -\frac{(t-a)^{\overline{\nu}}(b-s)^{\overline{\nu}}}{\Gamma(\nu+1)(b-a)^{\overline{\nu}}}+\frac{(t-s)^{\overline{\nu}}}{\Gamma(1+\nu)}\bigg] = -\frac{(t-a)^{\overline{\nu-1}}(b-s)^{\overline{\nu}}}{\Gamma(\nu)(b-a)^{\overline{\nu}}}+\frac{(t-s)^{\overline{\nu-1}}}{\Gamma(\nu)}.
	\]
This expression is nonnegative if and only if
	\[
	\frac{(t-a)^{\overline{\nu-1}}(b-s)^{\overline{\nu}}}{\Gamma(\nu)(b-a)^{\overline{\nu}}} \leq \frac{(t-s)^{\overline{\nu-1}}}{\Gamma(\nu)}.
	\]
Since $t-s$ is positive, this happens if
	\[
	\frac{(t-a)^{\overline{\nu-1}}(b-s)^{\overline{\nu}}}{(b-a)^{\overline{\nu}}(t-s)^{\overline{\nu-1}}} \leq 1.
	\]
Now, by definition of the rising function,
	\begin{align*}
	&\frac{(t-a)^{\overline{\nu-1}}(b-s)^{\overline{\nu}}}{(b-a)^{\overline{\nu}}(t-s)^{\overline{\nu-1}}}\\
	&= \frac{\Gamma(t-a+\nu-1)}{\Gamma(t-a)}\frac{\Gamma(b-s+\nu)}{\Gamma(b-s)}\frac{\Gamma(b-a)}{\Gamma(b-a+\nu)}\frac{\Gamma(t-s)}{\Gamma(t-s+\nu-1)}\\
	&=\frac{(t-s+\nu-1)(t-s+\nu)\cdots(t-a+\nu-2)}{(t-s)(t-s+1)\cdots(t-a-1)}\\
	&\quad \cdot \frac{(b-s)(b-s+1)\cdots(b-a-1)}{(b-s+\nu)(b-s+\nu+1)\cdots(b-a+\nu-1)}\\
	&=\frac{(t-s+\nu-1)}{(t-s)}\frac{(t-s+\nu)}{(t-s+1)}\cdots\frac{(t-a+\nu-2)}{(t-a-1)}\frac{(b-s)}{(b-s+\nu)}\frac{(b-s+1)}{(b-s+\nu+1)}\cdots\\
	& \quad \cdot \frac{(b-a-1)}{(b-a+\nu-1)}\\
	&\leq (1)(1)\cdots(1)(1)(1)\cdots(1)=1.
	\end{align*}
Next, we check $v(t,s)$ at the right endpoint, $t=b$, 
	\[
	v(b,s)= -\frac{(b-a)^{\overline{\nu}}(b-s)^{\overline{\nu}}}{\Gamma(\nu+1)(b-a)^{\overline{\nu}}} +\frac{(b-s)^{\overline{\nu}}}{\Gamma(\nu+1)} = 0.
	\]
Thus, $v(t,s)$ is nonpositive for $a\le s < t\le b$. Also note that $v(t,s)=u(t,s)$ for $t=s$. Therefore, for $t\in\N_a^b$, $G(t,s)$ is nonpositive. 

(2) Since we know that $v(t,s)$ is always increasing for $a\le s\le t\le b$ and that for $s=t$, $v(t,s)=u(t,s)$, it suffices to show that 
	\[
	u(t,s) \geq -\bigg(\frac{b-a}{4}\bigg)\bigg(\frac{\Gamma(b-a+1)}{\Gamma(\nu+1)\Gamma(b-a+\nu)}\bigg).
	\]
Let $a\le t\le s\le b$.  Then 
	\[
	G(t,s)=u(t,s)=-\frac{(t-a)^{\overline{\nu}}(b-s)^{\overline{\nu}}}{\Gamma(\nu+1)(b-a)^{\overline{\nu}}}\geq -\frac{(s-a)^{\overline{\nu}}(b-s)^{\overline{\nu}}}{\Gamma(\nu+1)(b-a)^{\overline{\nu}}}.
	\]
Note that for $\alpha\in\N_1$ and $0<\nu<1$, 
	\[
	\alpha^{\overline{\nu}}=\frac{\Gamma(\alpha+\nu)}{\Gamma(\alpha)}\leq \frac{\Gamma(\alpha+1)}{\Gamma(\alpha)}=\alpha^{\overline{1}}.
	\]
So
	\begin{align*}
	-\frac{(s-a)^{\overline{\nu}}(b-s)^{\overline{\nu}}}{\Gamma(\nu+1)(b-a)^{\overline{\nu}}} &\geq -\frac{(s-a)^{\overline{1}}(b-s)^{\overline{1}}}{\Gamma(\nu+1)(b-a)^{\overline{\nu}}}\\
	&\geq -\frac{(\frac{a+b}{2}-a)(b-\frac{a+b}{2})}{\Gamma(\nu+1)(b-a)^{\overline{\nu}}}\\
	&= -\frac{(b-a)(b-a)\Gamma(b-a)}{4\Gamma(\nu+1)\Gamma(b-a+\nu)}\\
	&= -\frac{(b-a)\Gamma(b-a+1)}{4\Gamma(\nu+1)\Gamma(b-a+\nu)}\\
	&= -\bigg(\frac{b-a}{4}\bigg)\bigg(\frac{\Gamma(b-a+1)}{\Gamma(\nu+1)\Gamma(b-a+\nu)}\bigg).
	\end{align*}
Therefore, 
	\[
	G(t,s)\geq -\bigg(\frac{b-a}{4}\bigg)\bigg(\frac{\Gamma(b-a+1)}{\Gamma(\nu+1)\Gamma(b-a+\nu)}\bigg).
	\]
	
(3) Consider
	\begin{align*}
	\int_a^b |G(t,s)|\nabla s &= \int_a^t |v(t,s)|\nabla s+\int_t^b |u(t,s)|\nabla s \\
	&= \int_a^t \bigg| -\frac{(t-a)^{\overline{\nu}}(b-s)^{\overline{\nu}}}{\Gamma(\nu+1)(b-a)^{\overline{\nu}}}+\frac{(t-s)^{\overline{\nu}}}{\Gamma(1+\nu)} \bigg|\nabla s\\
	&\quad + \int_t^b \frac{(t-a)^{\overline{\nu}}(b-s)^{\overline{\nu}}}{\Gamma(\nu+1)(b-a)^{\overline{\nu}}} \nabla s\\
	&= \int_a^t -\bigg[ -\frac{(t-a)^{\overline{\nu}}(b-s)^{\overline{\nu}}}{\Gamma(\nu+1)(b-a)^{\overline{\nu}}}+\frac{(t-s)^{\overline{\nu}}}{\Gamma(1+\nu)} \bigg]\nabla s \\
	&\quad + \int_t^b \frac{(t-a)^{\overline{\nu}}(b-s)^{\overline{\nu}}}{\Gamma(\nu+1)(b-a)^{\overline{\nu}}} \nabla s\\
	&=\int_a^b \frac{(t-a)^{\overline{\nu}}(b-s)^{\overline{\nu}}}{\Gamma(\nu+1)(b-a)^{\overline{\nu}}} \nabla s -\int_a^t \frac{(t-s)^{\overline{\nu}}}{\Gamma(\nu+1)}\nabla s\\
	&= -\frac{(t-a)^{\overline{\nu}}(b-s-1)^{\overline{\nu+1}}}{\Gamma(\nu+2)(b-a)^{\overline{\nu}}}\bigg|_{s=a}^{s=b} +\frac{(t-s-1)^{\overline{\nu+1}}}{\Gamma(\nu+2)}\bigg|_{s=a}^{s=t}\\
	&= \frac{(t-a)^{\overline{\nu}}(b-a-1)^{\overline{\nu+1}}}{\Gamma(\nu+2)(b-a)^{\overline{\nu}}} - \frac{(t-a-1)^{\overline{\nu+1}}}{\Gamma(\nu+2)}\\
	&= \frac{(t-a)^{\overline{\nu}}(b-a-1)(b-a)^{\overline{\nu}}}{\Gamma(\nu+2)(b-a)^{\overline{\nu}}} - \frac{(t-a-1)(t-a)^{\overline{\nu}}}{\Gamma(\nu+2)}\\
	&= \frac{(t-a)^{\overline{\nu}}}{\Gamma(\nu+2)}[b-a-1-(t-a-1)]\\
	&= \frac{(t-a)^{\overline{\nu}}(b-t)}{\Gamma(\nu+2)}\\
	&\leq \frac{(t-a)(b-t)}{\Gamma(\nu+2)}\\
	&\leq \frac{(\frac{a+b}{2}-a)(b-\frac{a+b}{2})}{\Gamma(\nu+2)}\\
	&= \frac{(b-a)(b-a)}{4\Gamma(\nu+2)}\\
	&= \frac{(b-a)^2}{4\Gamma(\nu+2)}.
	\end{align*}

(4) 
We can assume that $b-a>1$ as if $b-a=1$, it would be an initial value problem and not a boundary value problem. Taking the difference of $u(t,s)$ with respect to $t$, we have 
\[\nabla_t u(t,s) = 
\nabla_t \frac{-(t-a)^{\overline{\nu}}(b-s)^{\overline{\nu}}}{\Gamma(\nu+1)(b-a)^{\overline{\nu}}}
=\frac{-\nu(t-a)^{\overline{\nu-1}}(b-s)^{\overline{\nu}}}{\Gamma(\nu+1)(b-a)^{\overline{\nu}}}.\]
Since all the factors in the above equation are individually positive, we get that $\nabla u(t,s)$ is nonpositive. Let $t\in\mathbb{N}_a$. 
We therefore know 
\begin{align*}
&\int_a^b\left|\nabla_tG(t,s)\right|\nabla s \\
&= \int_a^{t-1}\left|\nabla_tG(t,s)\right|\nabla s + \int_{t-1}^b\left|\nabla_tG(t,s)\right|\nabla s\\
&=\int_a^{t-1}\left|\nabla_tv(t,s)\right|\nabla s + \int_{t-1}^b\left|\nabla_tu(t,s)\right|\nabla s\\
&=\int_a^{t-1}\left[\nabla_t\frac{-(t-a)^{\overline{\nu}}(b-s)^{\overline{\nu}}}{\Gamma(\nu+1)(b-a)^{\overline{\nu}}} + \nabla_t\frac{(t-s)^{\overline{\nu}}}{\Gamma(\nu+1)}\right]\nabla s\\
&\quad+ \int_{t-1}^b\nabla_t\frac{(t-a)^{\overline{\nu}}(b-s)^{\overline{\nu}}}{\Gamma(\nu+1)(b-a)^{\overline{\nu}}}\nabla s\\
&=\int_a^{t-1}\nabla_t\frac{-(t-a)^{\overline{\nu}}(b-s)^{\overline{\nu}}}{\Gamma(\nu+1)(b-a)^{\overline{\nu}}}\nabla s +\int_a^{t-1} \nabla_t\frac{(t-s)^{\overline{\nu}}}{\Gamma(\nu+1)}\nabla s \\
&\quad+ \int_{t-1}^b\nabla_t\frac{(t-a)^{\overline{\nu}}(b-s)^{\overline{\nu}}}{\Gamma(\nu+1)(b-a)^{\overline{\nu}}}\nabla s\\
&=\int_a^{t-1}\frac{-\nu(t-a)^{\overline{\nu-1}}(b-s)^{\overline{\nu}}}{\Gamma(\nu+1)(b-a)^{\overline{\nu}}}\nabla s +\int_a^{t-1}\frac{\nu(t-s)^{\overline{\nu-1}}}{\Gamma(\nu+1)}\nabla s \\
&\quad + \int_{t-1}^b\frac{\nu(t-a)^{\overline{\nu-1}}(b-s)^{\overline{\nu}}}{\Gamma(\nu+1)(b-a)^{\overline{\nu}}}\nabla s\\
&=\frac{-\nu(t-a)^{\overline{\nu-1}}}{\Gamma(\nu+1)(b-a)^{\overline{\nu}}}\left[\frac{-1}{\nu+1}(b-s-1)^{\overline{\nu+1}}\right]_{s=a}^{s=t-1}\\
&\quad+\frac{\nu}{\Gamma(\nu+1)}\left[\frac{-1}{\nu}(t-s-1)^{\overline{\nu}}\right]_{s=a}^{s=t-1}\\
&\quad +\frac{\nu(t-a)^{\overline{\nu-1}}}{\Gamma(\nu+1)(b-a)^{\overline{\nu}}}\left[\frac{-1}{\nu+1}(b-s-1)^{\overline{\nu+1}}\right]_{s=t-1}^{s=b}\\
&=\frac{\nu(t-a)^{\overline{\nu-1}}}{\Gamma(\nu+2)(b-a)^{\overline{\nu}}}\left[(b-t)^{\overline{\nu+1}}-(b-a-1)^{\overline{\nu+1}}\right]\\
&\quad -\frac{1}{\Gamma(\nu+1)}\left[(t-t+1-1)^{\overline{\nu}}-(t-a-1)^{\overline{\nu}}\right]\\
&\quad +\frac{-\nu(t-a)^{\overline{\nu-1}}}{\Gamma(\nu+2)(b-a)^{\overline{\nu}}}\left[(b-b-1)^{\overline{\nu+1}}-(b-t)^{\overline{\nu+1}}\right]\\
&= \frac{2\nu(t-a)^{\overline{\nu-1}}(b-t)^{\overline{\nu+1}}}{\Gamma (\nu+2) (b-a)^{\overline{\nu}}} +\frac{(t-a-1)^{\overline{\nu}}}{\Gamma (\nu+1)}-\frac{\nu(t-a)^{\overline{\nu-1}}(b-a-1)}{\Gamma (\nu+2)}.
\end{align*}
Suppose $t=b$. This would imply that
\begin{align*}\int_a^b\left|\nabla_tG(t,s)\right|\nabla s &= \frac{2\nu(b-a)^{\overline{\nu-1}}(0)^{\overline{\nu+1}}}{\Gamma (\nu+2) (b-a)^{\overline{\nu}}} +\frac{(b-a-1)^{\overline{\nu}}}{\Gamma (\nu+1)}\\
&\quad-\frac{\nu(b-a)^{\overline{\nu-1}}(b-a-1)}{\Gamma(\nu+2)}\\
&=\frac{(\nu+1)(b-a-1)^{\overline{\nu}}}{\Gamma(\nu+2)}-\frac{\nu(b-a)^{\overline{\nu-1}}(b-a-1)}{\Gamma(\nu+2)}.
\end{align*}
For $t=b$ and $b-a=2$, this becomes
\begin{align*}\int_a^b\left|\nabla_tG(t,s)\right|\nabla s &= \frac{(\nu+1)(1)^{\overline{\nu}}}{\Gamma(\nu+2)}-\frac{\nu(2)^{\overline{\nu-1}}(1)}{\Gamma(\nu+2)}\\
&=\frac{(\nu+1)\Gamma(\nu+1)}{\Gamma(\nu+2)}-\frac{\nu\Gamma(\nu+1)}{\Gamma(\nu+2)}\\
&=1-\frac{\nu}{\nu+1}\\
&=\frac{1}{\nu+1}\\
&\le\frac{2}{\nu+1}\\
&=\frac{b-a}{\nu+1}.
\end{align*}
For $t=b$ and $b-a=3$, we have
\begin{align*}\int_a^b\left|\nabla_tG(t,s)\right|\nabla s &= \frac{(\nu+1)(2^{\overline{\nu}})}{\Gamma(\nu+2)}-\frac{\nu(3^{\overline{\nu-1}})(2)}{\Gamma(\nu+2)}\\
&=\frac{(\nu+1)\Gamma(\nu+2)}{\Gamma(\nu+2)}-\frac{2\nu\Gamma(\nu+2)}{\Gamma(\nu+2)\Gamma(3)}\\
&=\nu+1-\nu\\
&=1\\
&=\frac{3}{3}\\
&\le\frac{b-a}{\nu+1}.
\end{align*}
For $t=b$ and $b-a\ge4$, the result holds since
\begin{align*}\int_a^b\left|\nabla_tG(t,s)\right|\nabla s &= \frac{(\nu+1)(b-a-1)^{\overline{\nu}}}{\Gamma(\nu+2)}-\frac{\nu(b-a)^{\overline{\nu-1}}(b-a-1)}{\Gamma(\nu+2)}\\
&=\frac{(\nu+1)(b-a-2+\nu)\cdots(2+\nu)}{\Gamma(b-a-1)}\\
&\quad-\frac{\nu\Gamma(b-a-1+\nu)(b-a-1)}{\Gamma(2+\nu)\Gamma(b-a)}\\
&=\frac{(\nu+1)(b-a-2+\nu)\cdots(2+\nu)}{\Gamma(b-a-1)}\\
&\quad-\frac{\nu\Gamma(b-a-1+\nu)}{\Gamma(2+\nu)\Gamma(b-a-1)}\\
&=\frac{(\nu+1)(b-a-2+\nu)\cdots(2+\nu)}{(b-a-2)!}\\
&\quad-\frac{\nu(b-a-2+\nu)\cdots(2+\nu)}{(b-a-2)!}\\
&=\frac{(b-a-2+\nu)\cdots(2+\nu)}{(b-a-2)!}\\
&=\frac{(b-a-1)(b-a-2+\nu)\cdots(2+\nu)}{(b-a-1)!}\\
&\le\frac{(b-a-1)(b-a-1)\cdots(3)}{(b-a-1)!}\\
&=\frac{\frac{1}{2}(b-a-1)!(b-a-1)}{(b-a-1)!}\\
&=\frac{b-a-1}{2}\\
&\le\frac{b-a}{\nu+1}.
\end{align*}
So the result holds if $t=b$ generally. Now, assume $t<b$. If $t=a+1$, then we have
\begin{align*}
\int_a^b\left|\nabla_tG(t,s)\right|\nabla s &= \frac{2\nu(1^{\overline{\nu-1}})(b-a-1)^{\overline{\nu+1}}+ (\nu+1)(0^{\overline{\nu}})(b-a)^{\overline{\nu}}}{\Gamma(\nu+2)(b-a)^{\overline{\nu}}}\\
&\quad
-\frac{\nu(1^{\overline{\nu-1}})(b-a-1)^{\overline{\nu+1}}}{\Gamma(\nu+2)(b-a)^{\overline{\nu}}}\\
&= \frac{2\nu\Gamma(\nu)(b-a-1)(b-a)^{\overline{\nu}}-\nu\Gamma(\nu)(b-a-1)(b-a)^{\overline{\nu}}}{\Gamma(\nu+2)(b-a)^{\overline{\nu}}}\\
&= \frac{2\nu\Gamma(\nu)(b-a-1)-\nu\Gamma(\nu)(b-a-1)}{\Gamma(\nu+2)}\\
&= \frac{\Gamma(\nu+1)(b-a-1)}{\Gamma(\nu+2)}\\
&= \frac{\Gamma(\nu+1)(b-a-1)}{(\nu+1)\Gamma(\nu+1)}\\
&= \frac{b-a-1}{\nu+1}\\
&\le \frac{b-a}{\nu+1}.
\end{align*}
If $t=a+2$, then 
\begin{align*}
&\int_a^b\left|\nabla_tG(t,s)\right|\nabla s \\
&= \frac{2\nu(2^{\overline{\nu-1}})(b-a-2)^{\overline{\nu+1}}+ (\nu+1)(1^{\overline{\nu}})(b-a)^{\overline{\nu}}-\nu(2^{\overline{\nu-1}})(b-a-1)^{\overline{\nu+1}}}{\Gamma(\nu+2)(b-a)^{\overline{\nu}}}\\
&= \frac{2\nu\Gamma(\nu+1)(b-a-2)^{\overline{\nu+1}}}{\Gamma(\nu+2) (b-a)^{\overline{\nu}}} + \frac{(\nu+1)(\nu)\Gamma(\nu)(b-a)^{\overline{\nu}}}{\Gamma(\nu+2)(b-a)^{\overline{\nu}}} \\
&\quad- \frac{\nu\Gamma(\nu+1)(b-a-1)^{\overline{\nu+1}}}{\Gamma(\nu+2) (b-a)^{\overline{\nu}}}\\
&= \frac{2\nu\Gamma(\nu+1)(b-a-1)(b-a-2)}{(\nu+1)(b-a-1+\nu)} + 1 - \frac{\nu(b-a-1)}{\nu+1}\\
&\le \frac{2\nu(b-a-2)}{\nu+1} + 1 - \frac{\nu(b-a-1)}{\nu+1}\\
&= \frac{2\nu(b-a-2)}{\nu+1} + \frac{\nu+1}{\nu+1} - \frac{\nu(b-a-1)}{\nu+1}\\
&= \frac{\nu(b-a-2)+1}{\nu+1}\\
&\le \frac{b-a-1}{\nu+1}\\
&\le \frac{b-a}{\nu+1}.
\end{align*}
If $t=a+3$, then
\begin{align*}
&\int_a^b\left|\nabla_tG(t,s)\right|\nabla s \\
&= \frac{2\nu(3^{\overline{\nu-1}})(b-a-3)^{\overline{\nu+1}}+ (\nu+1)(2^{\overline{\nu}})(b-a)^{\overline{\nu}}-\nu(3^{\overline{\nu-1}})(b-a-1)^{\overline{\nu+1}}}{\Gamma(\nu+2)(b-a)^{\overline{\nu}}}\\
&= \frac{2\nu\Gamma(2+\nu)\Gamma(b-a-2+\nu)\Gamma(b-a)}{\Gamma(3)\Gamma(b-a+\nu)\Gamma(\nu+2)\Gamma(b-a-3)} + (\nu+1)\\&\quad - \frac{\nu\Gamma(\nu+2)(b-a-1)}{\Gamma(\nu+2)\Gamma(3)}\\
&= \frac{\nu(b-a-1)(b-a-2)(b-a-3)}{(b-a-1+\nu)(b-a-2+\nu)} + (\nu+1) - \frac{\nu(b-a-1)}{2}\\
&\le \nu(b-a-3) + \nu+1 - \frac{\nu(b-a-1)}{2}\\
&= \frac{2\nu b - 2\nu a - 6\nu + 2\nu + 2 - \nu b + \nu a + \nu}{2}\\
&= \frac{\nu(b-a-3) + 2}{2}\\
&\le \frac{b-a-3 + 2}{2}\\
&= \frac{b-a-1}{2}\\
&\le \frac{b-a}{\nu+1}.
\end{align*}
Now suppose that $t = a+k$, where $k\in\mathbb{N}_4^{b-a-1}$. Then
\begin{align*}
&\int_a^b\left|\nabla_tG(t,s)\right|\nabla s \\
&= \frac{2\nu(k)^{\overline{\nu-1}}(b-a-k)^{\overline{\nu+1}}}{(b-a)^{\overline{\nu}}\Gamma(\nu+2)} + \frac{(\nu+1)(k-1)^{\overline{\nu}}}{\Gamma(\nu+2)} - \frac{\nu(k)^{\overline{\nu-1}}(b-a-1)}{\Gamma(\nu+2)}\\
&= \frac{2\nu\Gamma(k+\nu-1)\Gamma(b-a-k+\nu+1)\Gamma(b-a)}{\Gamma(k)\Gamma(b-a+\nu)\Gamma(\nu+2)\Gamma(b-a-k)} + \frac{(\nu+1)\Gamma(k-1+\nu)}{\Gamma(\nu+2)\Gamma(k-1)}\\&\quad - \frac{\nu\Gamma(k+\nu-1)(b-a-1)}{\Gamma(\nu+2)\Gamma(k)}\\
&= \frac{2\nu(\nu+2)\dots(\nu+k-2)(b-a-1)\dots(b-a-k)}{(k-1)! (b-a-1+\nu)\dots(b-a-(k-1)+\nu)} + \frac{(\nu+1)\dots(\nu+k-2)}{(k-2)!}\\
&\quad - \frac{\nu(\nu+2)\dots(\nu+k-2)(b-a-1)}{(k-1)!}\\
&\le \frac{2\nu(\nu+2)\dots(\nu+k-2)(b-a-k)}{(k-1)!} + \frac{(k-1)(\nu+1)\dots(\nu+k-2)}{(k-1)!} \\&\quad- \frac{\nu(\nu+2)\dots(\nu+k-2)(b-a-1)}{(k-1)!}\\
&= \frac{\nu(\nu+2)\dots(\nu+k-2)(2b-2a-2k-b+a+1)}{(k-1)!} \\
&\quad+ \frac{(k-1)(\nu+1)\dots(\nu+k-2)}{(k-1)!} \\
&= \frac{\nu(\nu+2)\dots(\nu+k-2)(b-a+1-2k) + (k-1)(\nu+1)\dots(\nu+k-2)}{(k-1)!}\\
&\le \frac{(1)(3)(4)\dots(k-1)(b-a+1-2k) + (k-1)(2)(3)\dots(k-1)}{(k-1)!}\\
&= \frac{\frac{1}{2}(k-1)!(b-a+1-2k) + (k-1)(k-1)!}{(k-1)!}\\
&= \frac{(b-a+1-2k) + 2(k-1)}{2}\\
&= \frac{b-a-1}{2}\\
&\le \frac{b-a}{\nu+1}.
\end{align*}
}

These properties of the Green's function are important in proving the existence and uniqueness of solutions of BVPs for nonlinear difference equations.

\end{document}